\documentclass[aos,MSNbibl,dvips]{arximspdf}
\usepackage{algorithm}
\usepackage{algorithmic}

\doi{10.1214/13-AOS1141} 
\volume{41}
\issue{4}
\pubyear{2013}
\firstpage{2149}
\lastpage{2175}

\makeatletter

\renewcommand{\epsilon}{\varepsilon}

\newtheorem{proposition}{Proposition}
\newtheorem{theorem}[proposition]{Theorem}
\newtheorem{lemma}[proposition]{Lemma}
\newtheorem{corollary}[proposition]{Corollary}

\newproclaim{remark}{Remark}
\newproclaim{definition}[proposition]{Definition}
\newproclaim{example}[proposition]{Example}

\newcommand{\err}{\operatorname{err}}
\newcommand{\loss}{\operatorname{loss}}
\newcommand{\vol}{\operatorname{vol}}
\newcommand{\prob}{\operatorname{Prob}}
\newcommand{\mean}{\mathbb{E}}
\newcommand{\spn}{\operatorname{span}}
\newcommand{\dist}{\operatorname{dist}}
\newcommand{\reals}{\mathbb{R}}
\newcommand{\argmin}{\operatorname{\arg\min}}
\newcommand{\argmax}{\mathop{\arg\max}}
\newcommand{\E}{\mathcal{E}}
\newcommand{\G}{\mathcal{G}}
\renewcommand{\P}{\mathcal{P}}

\makeatother

\begin{document}
\begin{frontmatter}

\title{Nearly optimal minimax estimator for high-dimensional sparse
linear regression}
\runtitle{Nearly optimal minimax estimator}

\begin{aug}
\author[A]{\fnms{Li} \snm{Zhang}\corref{}\ead[label=e1]{lzha@microsoft.com}}
\runauthor{L. Zhang}
\affiliation{Microsoft Research Silicon Valley}
\address[A]{Microsoft Research\\
1065 La Avenida\\
Mountain View, California 94043\\
USA\\
\printead{e1}} 
\end{aug}

\received{\smonth{1} \syear{2013}}
\revised{\smonth{6} \syear{2013}}

%
\begin{abstract}
We present estimators for a well studied statistical estimation
problem: the estimation for the linear regression model with soft
sparsity constraints ($\ell_q$ constraint with $0<q\leq1$) in the
high-dimensional setting. We first present a family of estimators,
called \emph{the projected nearest neighbor estimator} and show, by
using results from Convex Geometry, that such estimator is within a
logarithmic factor of the optimal for any design matrix. Then by
utilizing a semi-definite programming relaxation technique developed in
[\textit{SIAM J. Comput.} \textbf{36} (2007) 1764--1776], we obtain an
approximation algorithm for computing the minimax risk for any such
estimation task and also a polynomial time nearly optimal estimator for
the important case of $\ell_1$ sparsity constraint. Such results were
only known before for special cases, despite decades of studies on this
problem. We also extend the method to the adaptive case when the
parameter radius is unknown.
\end{abstract}

%

\end{frontmatter}

\section{Introduction}

In the classical estimation problem with linear regression model, one
observes a noisy $\widetilde{y}$ of some $y\in\reals^n$
where $y=X\theta$ for a given $n\times p$ matrix $X$ (called the
\emph{design matrix}) and an unknown $\theta\in\reals^p$ and wishes
to estimate $y$ or $\theta$. Recently, there have been enormous
interests in the high-dimensional setting which in addition assumes
that the design matrix is high-dimensional, that is, when $p\gg n$, and
$\theta$ satisfies certain sparsity constraints. Such sparsity
constraints can be ``hard'', when it bounds the number of nonzero
components in $\theta$, or ``soft'', when $\theta$ is assumed to
belong to the unit $\ell_q$ ball for $0<q\leq1$. In the existing
study, the focus has so far been on the condition needed for $X$ such
that certain (typically polynomial time) estimators are nearly optimal
or achieve lowest possible error for the given parameters. The work
along this line has been quite
successful \cite
{dj-94,baraud-00,bbm-99,b-04,d-04,ct-05,d-06,ct-07,brt-09,cai-10,bm-10,z-10,
zhangtong-11,rwy-11,djmm-11}\vadjust{\goodbreak}
and produced many characterization of $X$ (typically Gaussian random
matrix) for which a polynomial time nearly optimal estimator exists.

The main departure point of this study is that we consider the problem
of designing nearly optimal estimator for \emph{any} given design
matrix $X$, that is, we make no assumption about $X$. As the main
contribution of this paper, we present a family of estimators, which
we call \emph{the projected nearest neighbor estimator} (PNN), and show
that for any design matrix $X$, there is a projected
nearest neighbor estimator that is nearly optimal in terms of the
prediction risk for the corresponding linear regression problem over
soft sparsity constraints. As a consequence, we obtain a
polynomial time algorithm to compute the approximate minimax risk for any
such problem and a polynomial time estimator in the important case of
$q=1$. Our results represent the first provably nearly optimal
estimators without any constraint on the design matrix for $0<q\leq1$.
We also design an adaptive estimator for the case when the $\ell_1$ radius
is not given.

We believe that studying optimal estimator for arbitrary $X$ is
important for multiple reasons. First, in practice we often do not
have control over the design matrix or even the distribution of the
design matrix. The design matrix might be ``ill''-conditioned such
that no estimator can achieve good accuracy. On the other hand, the
design matrix may have a structure, as is often the case in practice,
rather than completely random. In this case, it is important to take
advantage of such structure to obtain better accuracy. Second,
while there have been many characterization (typically some isometry
property on $X$) known for certain algorithms to work well, it is
often difficult to tell if the required property holds for a given
$X$. So most results assume that $X$ come from Gaussian random matrix.
Third, relaxing the requirement about the design $X$ calls for the
development of new algorithms as well as new analysis tools. Indeed, to
argue the optimality of our estimator, we have to utilize novel tools
from Convex Geometry (the classical
restricted invertibility result by Bourgain and
Tzafriri \cite{bt-87}).

\subsection{Problem setup}

In the linear regression problem, one observes
$\widetilde{y} = y+g\in\reals^n$, where $y=X\theta$ for a given
$n\times p$ matrix $X$ and an unknown vector $\theta\in\ell_q(C)$
for $0<q\leq1$, where
$\ell_q(C)=\{(\theta_1,\ldots,\theta_p)\dvtx  (\sum_i|\theta
_i|^{q})^{1/q}\leq C\}$. In
addition, the noise $g$ is a random vector drawn from the multivariate
Gaussian distribution with the covariance matrix $\sigma^2I$. In this
paper, we only consider the prediction estimation, that is, on the
estimation of $y$ but not $\theta$. We use the standard total squared
loss\footnote{We use the total squared error instead of the common
mean squared error purely for the brevity of notation.} to measure the
error of an estimation, that is,
\[
\loss(\widehat{y}, y) =\|\widehat{y}-y\|^2=\sum
_i(\widehat{y}_i - y_i)^2.
\]

For an estimator $M\dvtx \reals^n\to\reals^n$, we define the expected
error of $M$ on an input $y$ and on Gaussian error as
\[
\err_M(y,\sigma) = \mean_{\widetilde{y}=y+g; g\sim\G(\sigma
)}\loss\bigl(M(\widetilde{y}),
y\bigr) = \mean_{\widetilde{y}=y+g; g\sim\G
(\sigma)}\bigl\|M(\widetilde{y})-y\bigr\|^2.
\]

Following \cite{dlm-90}, for $K\subseteq\reals^n$, the risk of $M$
over $K$ is defined as
%
\begin{equation}
\label{eqR} R_M(K,\sigma) = \sup_{y\in K}
\err_M(y,\sigma).
\end{equation}

Define the minimax risk, denote by $R^\ast(K,\sigma)$, as the minimum
achieve-able risk among all the possible estimators, that is,
%
\begin{equation}
\label{eqRs} R^\ast(K,\sigma) = \inf_M
R_M(K,\sigma).
\end{equation}

For the aforementioned linear model with sparsity constraint
$\ell_q(C)$, we have $K=X\ell_q(C)$ for an $n\times p$ design
matrix $X$. Clearly, the minimax risk $R^\ast$ ranges between $0$ and
$n\sigma^2$ and depends on the structure of $X$. The main goal of
this paper is to design an estimator $M$ such that
$R_M(X\ell_q(C),\sigma)$ is close to $R^\ast(X\ell_q(C),\sigma)$ for
any given $X$. For our main results, we consider the case where the
sparsity radius $C$ is given. Since we will only consider the
prediction risk, we can assume, by rescaling $X$, that $C=1$. In what
follows, we write $\ell_q$ for $\ell_q(1)$. In addition, we only
consider the high-dimensional case where $p\geq n$ because for
$p<n$, we can apply a rotation to the design matrix so that the last
$n-p$ rows are entirely $0$. Since Gaussian noise is invariant
under rotation, this does not affect the minimax risk, and the
dimensions of the design matrix is effectively reduced to
$p\times p$.

\subsection{Main contribution}

We present a family of estimators, called the \emph{projected nearest
neighbor estimator} (PNN), that can achieve nearly optimal risk for
\emph{any} design matrix $X$ and any given $0<q\leq1$. The projected
nearest neighbor estimator is a combination of two classic estimators:
the \emph{orthogonal projection estimator}, in which the estimation is
obtained by projecting the observation $\widetilde{y}$ to a properly
chosen subspace, and the \emph{nearest neighbor estimator}, in which
$\widetilde{y}$ is mapped to the closest point (in terms of $\ell_2$
distance) on the ground truth set $K$. The projected nearest neighbor
estimator is defined with respect to an orthogonal projection $P$. It
is the summation of two components: one, similar to the orthogonal
projection estimator, is the projection $P\widetilde{y}$ of
$\widetilde{y}$ by $P$; the other, similar to the nearest neighbor
estimator, is the nearest neighbor projection of $P^\bot \widetilde{y}$
on $P^\bot K$, where $P^\bot$ is the projection orthogonal to $P$. As
the main contribution of this work, we show that for any $X$,
$0<q\leq1$, and $\sigma>0$, there always exists a projection $P$ so
that the corresponding projected nearest neighbor estimator for
$K=X\ell_q$ is nearly minimax optimal. More precisely, we show the
following theorem.\footnote{Throughout this paper, the $O$ notation
only hides some absolute constant, that is, a constant independent of
any of the parameters, such as $n,p,q,\sigma,X,\theta,y$.}
%
\begin{theorem}\label{thmmain}
For any given $n\times p$
matrix $X$, $0<q\leq1$, and $\sigma\geq0$, there exists a projected
nearest neighbor estimator $M$ such that
\[
R_M(X\ell_q,\sigma) = O\bigl(c_{q} \bigl(
\log^{1-q/2}p\bigr) R^\ast(X\ell_q,\sigma)\bigr),
\]
where $c_{q} = O(2^{{1}/{q}}\frac{1}{q}\ln\frac{2}{q})$ is a
constant dependent on $q$ only.
\end{theorem}

In the above theorem, the projection $P$ is
chosen in two steps: (1) for each $0\leq k\leq n$, a $k$-dimensional
projection $P_k$ is chosen to minimize $\max_i \| P^{\bot} x_i \|$
where $x_i$'s are column vectors of $X$; (2) a proper $k^\ast$ is
chosen to minimize the risk among all the $P_k$'s. Finding the
projection in step 1 turns out to be NP-hard. However, by using the
semi-definite programming technique in \cite{vvyz-07}, we can compute
an approximately optimal projection and therefore an approximate
minimax risk in polynomial time.
%
\begin{theorem}\label{thmalgo}
For any given $n\times p$
matrix $X$, $0<q\leq1$, and $\sigma\geq0$, we can compute an
$O(c_{q} \log p)$
approximation\footnote{A quantity $a$ is a $c$-approximation of
$a^\ast\geq0$, if $a^\ast\leq a \leq c a^\ast$.} of $R^\ast(X\ell_q,
\sigma)$ in polynomial time. When $q=1$, there is a randomized
polynomial time estimator that is within $O(\log p)$ factor of the
optimal.
\end{theorem}

The above two results assume that the radius of $\ell_q$ ball is
given. For $q=1$, we can extend the estimator to the adaptive case
when $\|\theta\|_1$ is unknown. Using the similar idea to the
projected nearest neighbor estimator, we have that
%
\begin{theorem}\label{thmadaptive}
There is a polynomial time adaptive estimator $A$ such that for any
given $n\times p$ matrix $X$, $\theta$, and $\sigma>0$,
%
\begin{equation}
\label{eqada} \err_A(X\theta, \sigma) = O\bigl(\log p\cdot
R^\ast\bigl(X \ell_1\bigl(\|\theta\|_1\bigr), \sigma
\bigr)+\sqrt{n\log n} \sigma^2\bigr).
\end{equation}
\end{theorem}

Notice that the first term of the above error is $O(\log p)$ factor
within the oracle risk bound when $\|\theta\|_1$ is given. While we do
not quite get the true oracle bound due to the presences of the
additive term of $\sqrt{n\log n} \sigma^2$, the bound becomes a true
(and nontrivial) oracle bound for a rather large range of
$\|\theta\|_1$. See Remark \ref{rmkada} for a more detailed
discussion.

\subsection{Intuition}

We provide some high level intuition of the projected nearest neighbor
estimator. The orthogonal projection estimator, by projecting the
observation to a chosen subspace, effectively identifies the ``leading
factors'' in the ground truth set. It works well when $K$ is
``skewed''. However by simple projection, it ignores the detailed
local geometry of $K$. This makes it less effective when $K$ has many
constraints or has constraints involving many dimensions, for example, when
$K$ satisfies sparse constraints. On the other hand, the nearest
neighbor estimator, by projecting to the nearest neighbor, depends
more on the local geometry of $K$. But it ignores the global geometry
of $K$ so it works well when the body is not skewed along any
direction. In some sense, the projected nearest neighbor estimator
achieves the optimality by taking both global and local geometry into
account: it first identifies the skewed dimensions and then applies
the nearest neighbor estimator to the ``residual'' space which is less
biased.

It is long known that the nearest neighbor estimator may be far away from
the optimal when there is strong correlation among column vectors
of the design matrix
$X$ \cite{fan-penlike-01,fan-penlike-04,zou-lasso-06}. There have
been many methods proposed to deal with this problem. The projection
phase can be viewed as one way to remove the correlation such that the
residual vectors are less biased. This might not be
obvious as the projection only minimizes the maximum of $\ell_2$ norm
of the projection, a seemingly different quantity. However, in order
for the projected vectors to be all short, they necessarily
``span'' all the directions because otherwise we could ``tilt'' the
projection to reduce the longest projection. This intuition can
actually be made rigorous with the help of tools from Convex
Geometry \cite{bt-87}.

The technical analysis of the projected nearest neighbor estimator is
inspired by two recent works, one is the analysis on the nearest neighbor
estimator by Raskutti, Wainwright and Yu \cite{rwy-11}; the other is
on the optimality of the orthogonal projection estimator by Javanmard and
the author \cite{jz-12}. In \cite{rwy-11}, it is shown that if $X$
satisfies a certain isometry property, then the nearest neighbor estimator
is close to optimal. On the other hand, \cite{jz-12} shows that for
symmetric convex bodies there always exists a
projection such that the orthogonal projection estimator is close to
optimal. At the very high level, we combine the analysis of these two
results and show that there always exists a nearly optimal projection
of $X$ such that the bound in \cite{rwy-11} is nearly optimal on the
projected body.

While the main machinery in our analysis is similar to what is
in \cite{rwy-11} and \cite{jz-12}, we need further insights for our
problem. For the nearest neighbor analysis, we need a slightly different
analysis than \cite{rwy-11} to obtain an upper bound suit our
purpose. This also allows our result hold for all ranges of $p,n$.
The lower bound is obtained by extending the techniques
in \cite{jz-12} to the sets of the form $X\ell_q$ for $0<q\leq1$.
The technique utilizes some classical results from Banach space
geometry, first started by Bourgain and Tzafriri \cite{bt-87} and
fully developed by Szarek, Talagrand, and
Giannopoulous \cite{st-89,g-95}.

Despite its somewhat involved analysis, the projected nearest neighbor
estimator suggests a quite natural heuristic: project $K=X\ell_q$
to a subspace to make it more ``round'' before applying other
estimators (in our case the nearest neighbor estimator). This approach
is probably already being used in practice.
As the main result in this paper, we prove that such heuristics can
actually lead to a nearly optimal estimator. In addition, a nearly
optimal projection can be found in polynomial time via semi-definite
programming technique in \cite{vvyz-07}.

For the adaptive estimator, we consider the case of $q=1$. The
well-known Lasso \cite{t-lasso-96} and Dantzig selector \cite{ct-07}
can be
viewed as the adaptive version of the nearest negibhor estimator.
According to \cite{brt-09}, these estimators can achieve an error
bound dependent on $\|\theta\|_1$, which is the same as the oracle
risk bound of PNN when the projection is taken as the identity
projection. We can apply Lasso or Dantzig selector to the projection
of $X$ and to obtain the oracle risk bound of PNN under different
projection dimensions. This way, we can obtain a set of estimations
among which one achieves the true oracle risk bound! Unfortunately, we
cannot reliably determine which one it is. By using ideas from
hypothesis testing, we can only choose one within
$O(\sqrt{n\log n}\sigma^2)$ error, which accounts for the additive
bound in Theorem \ref{thmadaptive}.

More concretely, in PNN, the optimal projection dimension is a
staircase function of the parameter radius. So we try to ``guess''
$\|\theta\|_1$ at those critical values at which the optimal dimension
changes value. The problem then reduces to a hypothesis testing
problem on whether $y=X\theta$ belongs to some convex body. By using
the statistics of $\|\widetilde{y}-\widehat{y}\|_2^2$, we can achieve
the claimed bound. Our procedure is similar in spirit to the
classical Lepski's recipe \cite{l-adaptive-91,b-lepski-01} for
converting a nonadaptive estimator to an adaptive one. But there is a
significant difference as the PNN estimator is nonlinear, and the
projections at different dimensions lack a nested structure. As a
result, our bound leaves an additive gap of $\sqrt{n\log n}
\sigma^2$.

\subsection{Related work}

There are vast amounts of work on the minimax risk estimator. We refer
to \cite{nemirovski-98,tsybakov-09,johnstone-11} for comprehensive
surveys. Despite many studies on this subject, optimal or nearly
optimal estimators are only known for special types of bodies.

One particularly interesting case is when the parameter space is sparse.
It
is long known that no linear estimator works well under such
constraints (see, e.g., \cite{dlm-90}). Instead, one needs
nonlinear estimator such as the thresholding estimator to achieve
nearly optimal risk. Recently, much attention has been paid to the
(hard) sparsity constraint defined as the number of nonzero
components, dubbed as $\ell_0$ quantity, of a vector. This problem,
called
\emph{compressive sensing} in the literature, is computationally
infeasible in general so the study has focused on the condition under
which nearly optimal polynomial time estimator
exists \cite
{baraud-00,bbm-99,b-04,d-04,ct-05,ct-07,cai-10,bm-10,z-10,zhangtong-11}.

The case of $q=1$ is closely related to
Lasso \cite{t-lasso-96}, which is the nearest neighbor estimator for
the case of $q=1$ and later evolves to solving a regularized nearest
neighbor problem with the
$\ell_1$ norm penalty. While Lasso has proved to be very effective,
it is known that when the design matrix has strong correlation, the
Lasso estimator may not produce a good
estimation \cite{fan-penlike-01,fan-penlike-04}. Various methods have
been proposed to remove the
correlations \cite{fan-penlike-01,fan-penlike-04,zou-lasso-06} by
using different penalty terms. The projected nearest neighbor estimator
can also be viewed as a way to remove correlation. The difference is
that our method can be shown to be close to the optimal solution for
any design matrix $X$. In the projected nearest neighbor estimator, we
choose the projection dimension that balance two error terms. Similar
technique has appeared before. For example, in \cite{barron-greedy-08},
the estimation is chosen among greedy approximations of the span of
vectors of varying size, and the optimal choice is by balancing two error
terms. In \cite{cai-omp-11}, the dimension is controlled by a
stopping rule dependent on the noise structure. Despite these
similarity, the optimality of the projected nearest neighbor estimator
requires careful choice of the projection via solving a semi-definite
program. It is unlikely that the greedy algorithm can achieve the
same goal. On the other hand, the computational efficiency of the
greedy algorithm makes it (or some variation) an attractive practical
alternative to the more complex projection phase in this paper.

Many authors also consider (arguably more flexible and realistic) soft
sparsity constraints in the form of $\theta\in\ell_q$ for $0<q\leq
1$, the setting considered in this paper. In \cite{dj-94},
asymptotically tight bounds are obtained for $X=I$, the identity
matrix. A~similar notion of roughness was studied in \cite{lg-02} in
which soft-thresholding estimator is shown to be nearly optimal, again
for $X=I$, but extended to more general noise and loss
models. In \cite{d-06}, it is shown that there exists design matrices
$X$ which allow fairly accurate estimation when there is no noise.
In \cite{yz-lasso-10}, the authors presented several upper bounds,
dependent on the design matrix $X$, on the loss of the Lasso and
Dantzig selector methods when applied to soft sparsity constraints.
They also show that the upper bound is nearly optimal for a family of $X$'s.
Then in~\cite{rwy-11}, it is shown that the nearest neighbor estimator
is nearly optimal if $X$ satisfies certain isometry property which
holds for Gaussian random matrix $X$. In \cite{djmm-11}, it is shown
that for Gaussian random matrix, the (polynomial time) $\ell_1$ penalized
least squares is nearly optimal. Despite all these studies, no nearly
optimal estimator is known for general design matrix $X$. So our
knowledge is limited to the case where $X$ is a diagonal matrix or
when $X$ satisfies strong isometry properties. In \cite{cd-11}, the
authors showed a lower bound of the minimax risk on the estimation of
$\theta$ for any design matrix and with the hard sparsity constraint,
but it could be far away from the upper bound in general.

Among previous work, \cite{rwy-11} is particularly relevant to our
current work. In \cite{rwy-11}, the authors show, among many other
results, an upper bound for the nearest neighbor estimator which depends
on $q$ and the radius of $K$. While this could be far away from the
optimal, it turns out if we apply proper projection of $K$, the radius
of the projection can be made so that the resulted bound is near
optimal. For this, we follow a similar approach as \cite{jz-12}, in
which they show that the orthogonal projection estimator is nearly optimal
for symmetric linear constraints. But we need to adapt the argument
in \cite{jz-12} as $X \ell_q$ have exponentially many faces and can be
nonconvex.

As mentioned earlier, the transformation from nonadaptive estimator
to the adaptive one is similar to Lepski's
method \cite{l-adaptive-91,b-lepski-01} but there are significant
differences as our nonadaptive estimator does not quite satisfy the
properties required by Lepski's method.

\section{Preliminaries}

\subsection{Basic notation and definitions}

For a vector $x = (x_1,\ldots,x_p)\in\reals^p$ and $q>0$, denote by
$\|x\|_{q} = (\sum_i |x_i|^{q} )^{1/q}$. When $p\geq1$,
$\|x\|_{q}$ is a norm. When $0<q<1$, $\|x\|_{q}$ is not a norm but it is
quasi-convex as there is a constant $c$ dependent on $q$ such that
for any $x,y$, $\|x+y\|_{q}\leq c(\|x\|_{q}+\|y\|_{q})$. We use
$\ell_q^p(r)$ to denote the $p$-dimensional $q$-ball with radius
$r$, that is,
\[
\ell_q^p(r) = \bigl\{ x \in\reals^p\dvtx
\|x\|_{q} \leq r\bigr\}.
\]

We often drop $p$ when the dimension is clear from the context. We
use $\ell_q$ as a short hand for $\ell_q(1)$. For a set
$K\subseteq\reals^n$ containing the origin, define the $q$-radius of $K$
as $\|K\|_{q}=\sup_{x\in K} \|x\|_{q}$. In all these notations, whenever
$q$ is omitted, it means $q=2$.

We use $\G^n(\sigma)$ to denote the distribution of $n$-dimensional
Gaussian random variable with covariance matrix $\sigma^2 I$. Again,
we often drop $n$ and $\sigma$ when they are clear from the context.

As standard, $f = O(g)$ if there exists a constant $c>0$ such that
$f\leq c \cdot g$ and $f=\Omega(g)$ if there exists a constant $c>0$
such that $f\geq c \cdot g$. Throughout this paper, high probability
is understood as the probability of $1-1/n^2$.

\subsection{Minimax risk}\label{secprelimminimax}

An estimator $M$ is a map from $\reals^n$ to $\reals^n$: it takes a
noisy observation $\widetilde{y} = y+g$ of an unknown vector
$y\in\reals^n$ and maps it to an estimation $\widehat{y} =
M(\widetilde{y})$. Here we consider the noise drawn from
$\G^{n}(\sigma)$. As described early, the risk $R_M(K,\sigma)$ of
$M$ is defined as the maximum expected error among $y$ in $K$, that is,
\[
R_M(K,\sigma) = \sup_{y\in K} \mean_{\widetilde{y}=y+g; g\sim\G
(\sigma)}
\bigl[\bigl\|M(\widetilde{y})-y\bigr\|^2\bigr].
\]

The minimax risk of $K$ is defined as the minimum achievable risk for
$K$, that is, $R^\ast(K,\sigma) = \inf_{M} R_M(K,\sigma)$. We state
a well-known lower bound on the minimax risk of Euclidean balls which
we will
use later.
%
\begin{lemma}\label{lemlower}
$R^\ast(\ell_2^n(r), \sigma) = \Omega(\min(n\sigma^2, r^2))$.
\end{lemma}

\subsection{Orthogonal projection estimator}

The orthogonal projection estimator $T$ is a special type of linear
estimator. It is defined with respect to some linear subspace. The
estimation is simply by projecting the observation
$\widetilde{y}\in\reals^n$ to the subspace. Let $\P_k$ denotes all
the $k$-dimensional linear subspaces in $\reals^n$. For $P\in\P_k$,
we also use $P$ denote the orthogonal projection to $P$. The estimator
$T_P$ is then defined as $T_P(\widetilde{y}) = P\widetilde{y}$.

Since Gaussian random vector is invariant under the rotation, we have
that $R_{T_P}(K,\sigma)= k\sigma^2 + \sup_{y\in K} \|y-Py\|^2 =
k\sigma^2 + \sup_{y\in K} \|P^{\bot}y\|^2$, where $P^{\bot}$ denotes
the $(n-k)$-dimensional subspace orthogonal to $P$. For $0\leq k\leq
n$, define \emph{Kolmogorov width} (as in \cite{p-84}) as
\[
d_k(K)=\inf_{P\in\P_k} \sup_{y\in K}
\|y-Py\|.
\]

For $\ell_2$ norm, this definition is equivalent to following
more convenient form, which we will use through the paper:
\[
d_k(K)=\inf_{P\in\P_k}\bigl\|P^\bot(K)\bigr\|=\inf
_{P\in\P_{n-k}} \bigl\| P(K)\bigr\|.
\]

Clearly, $d_k(K)$ is monotonically decreasing with $k$. Kolmogorov
width determines the minimax risk of the orthogonal projection
estimators \cite{dlm-90}. Let $R_T$ denote the minimum risk among all
the orthogonal projection estimators.

\begin{lemma}\label{lemtrunx}
$R_T(K,\sigma) = \min_k (k\sigma^2 + d_k(K)^2)$.
\end{lemma}

The orthogonal projection estimator is long known to be nearly optimal
for ellipsoids \cite{pinsker-80,ih-81} and more generally for
quadratically convex and orthosymmetric objects \cite{dlm-90}.
However, it is also well known that the orthogonal projection
estimator (actually any linear estimator) can be far away from optimal
for the
$\ell_1$ ball and therefore does not work well for linear regression
with sparsity constraints.
%
\begin{lemma}[(\cite{dlm-90})]
\[
R_T\bigl(\ell_1^n, 1/\sqrt{n}\bigr) =
\Omega\bigl(\sqrt{n/\log n} R^\ast\bigl(\ell_1^n,1/
\sqrt{n}\bigr)\bigr).
\]
\end{lemma}

\subsection{Nearest neighbor estimator}

The nearest neighbor estimator is another well-known estimator. It
maps an observation to the nearest point on $K$,
that is, $N_K(\widetilde{y}) = \argmin_{\widehat{y}\in K} \|
\widehat{y} - \widetilde{y}\|$. The nearest neighbor estimator is
a nonlinear estimator and works well for ``skinny'' objects such as
the $\ell_1$ ball. However, we can construct an example
(Section \ref{subsecbadexample}) to demonstrate it is far from
optimal. Denote by $R_N(K,\sigma)$ the risk of the nearest neighbor
estimator.
%
\begin{lemma}\label{lemnnbad}
There exist ellipsoids $E_{n}\subset\reals^{n}$ for
$n=1,2,\ldots$ such that $R_N(E_{n}, 1) = \Omega(\sqrt{n}
R^\ast(E_{n},1))$.
\end{lemma}

\section{Projected nearest neighbor estimator}

We now describe the projected nearest neighbor estimator, which is
defined with respect to some low-dimensional orthogonal projection.
Given a $k$-dimensional subspace $P\in\P_k$, we define the projected
nearest neighbor estimator $H_P$ as follows. Let $P^{\bot}$ denote the
$(n-k)$-dimensional subspace orthogonal to $P$. Recall that we also
use $Px$, $P^\bot x$ to denote, respectively, the orthogonal
projection to the space $P$ and $P^\bot$. The estimator $H_P$ is
defined as
\[
H_P(\widetilde{y}) = P\widetilde{y} + N_{P^{\bot}K}
\bigl(P^{\bot
}\widetilde{y}\bigr).
\]

In other words, $H_P$ consists of two components, one of which is the
projection to the subspace $P$ and the other the nearest neighbor of
$P^{\bot}\widetilde{y}$ to $P^{\bot} K$. We use $R_H(K,\sigma)=\inf_{q}
R_{H_P}(K,\sigma)$ to denote the minimum risk achievable by the
projected nearest neighbor estimator for given $K,\sigma$.

When the projection is set as the identity projection, the
corresponding PNN is the same as the
nearest neighbor estimator. In addition, for the same projection, the
projected nearest
neighbor estimator outperforms the corresponding orthogonal projection
estimator. So the projected nearest neighbor estimator
subsumes both the nearest neighbor and the orthogonal projection
estimators. In the following, we give an example to show the projected
nearest neighbor
estimator can outperform both the orthogonal projection and the nearest
neighbor estimators by a large factor.

\begin{example}
Consider the ellipsoid defined as
\[
E_{n,k} = \Biggl\{ x\dvtx\frac{1}{\sqrt{n}}\sum
_{i=1}^k x_i^2 + \sum
_{i=k+1}^nx_i^2
\leq1\Biggr\}.
\]

Let $K = E_{n^2,k} \times\ell_1^{n}(\sqrt{n})$ with $k\leq n
$. By the
above discussion, we can see that for the orthogonal projection
estimator $R_T(K,1) = \Theta(n)$, and for the nearest neighbor
estimator $R_N(K,1)=\Theta(n)$, but $R_H(K,1) = O(\sqrt{n\log n
})$ by
setting $P$ to be the $k$-dimensional projection spanned by the $k$
long axes of $E_{n^2, k}$. This demonstrates a large gap between the
projected nearest neighbor estimator and both the orthogonal
projection and the nearest neighbor estimators.
\end{example}

To study the performance of the projected nearest neighbor estimator, we
first need the following error bound for the nearest neighbor estimator
from~\cite{rwy-11}.
%
\begin{proposition}\label{proplp}
For $0<q\leq1, K=X\ell_q^p$, the nearest neighbor estimator $N$ has
risk
\[
R_N(K,\sigma) = O\bigl(c_{q} \|K\|^{q}
\sigma^{2-q} (\log p)^{1-q/2}\bigr),
\]
where $c_q= O(2^{{1}/{q}}\frac{1}{q}\ln\frac{2}{q})$ is a
constant dependent on $q$ only.
\end{proposition}

The above bound is almost identical to Theorem 4(a)
in \cite{rwy-11}. We will present a slightly different proof which
applies to wider combination of parameters. For clarity and
completeness, we present the proof in Section \ref{subsecproplp}.
According to Proposition \ref{proplp}, the error is bounded by
$\|K\|^{q}$. Hence, if we fix the dimension of the projection in a
PNN estimator, in order to minimize the risk, we should seek the
projection $P$ that minimizes $\|PK\|$, that is, realizes Kolmogorov
width. By using this projection, we obtain the following upper bound
of the projected nearest neighbor estimator.
%
\begin{corollary}\label{corupperbound}
For any $0<q\leq1$ and any $K=X\ell_q^p$,
%
\begin{equation}
\label{equpperbound} R_H(K,\sigma) = O\Bigl(\min
_{0\leq k\leq n} \bigl(k\sigma^2 + c_{q}
d_k(K)^{q} \sigma^{2-q} (\log p)^{1-q/2}
\bigr)\Bigr),
\end{equation}
where $c_q$ is the same as in Proposition \ref{proplp}.
\end{corollary}
\begin{pf}
For any fixed $k$, the error consists of two terms: $O(k\sigma^2)$ for
the projection, and $O(c_{q} d_k(K)^{q} \sigma^{2-q} (\log p)^{1-q/2})$
for the nearest neighbor estimation. The second term comes from
Proposition \ref{proplp} with $\|K\|$ replaced by $d_k(K)$ if we apply
the projection that realizes $d_k(K)$. Clearly, we can choose $k$
with the minimum bound.
\end{pf}

To show (\ref{equpperbound}) is nearly optimal, we prove an almost
matching lower bound in terms of the Kolmogorov width. This is the key
technical contribution of the paper and relies on the classic
restricted invertibility property developed by Bourgain and
Tzafriri \cite{bt-87}. The proof is in Section \ref{subsecthmlowerbound}.
%
\begin{theorem}\label{thmlowerbound}
For $K=X\ell_q^p$,
%
\begin{equation}
\label{eqlowerbound} R^\ast(K,\sigma) = \Omega\Bigl(\max
_{0\leq k\leq n} \min\bigl(k\sigma^2, k^{1-2/q}
d_k(K)^2\bigr)\Bigr).
\end{equation}
\end{theorem}

Theorem \ref{thmmain} follows readily from
Corollary \ref{corupperbound} and Theorem \ref{thmlowerbound} by
setting $k$ to equalize two terms in (\ref{eqlowerbound}). The
details are in Section \ref{subsecthmmain}.

\begin{remark}
In the proof of Theorem \ref{thmmain}, we choose $k^\ast$ such that
$d_k(X)\approx k^{1/q} \sigma$. When $q$ goes to $0$, then $k^\ast$
goes to $1$. Therefore, when $q$ is close to $0$, the projected
nearest neighbor estimator becomes the ordinary nearest neighbor
algorithm. As stated in Theorem 4(b) in \cite{rwy-11}, the risk of
the nearest neighbor estimator is $O(s \log(p/s)\sigma^2 )$ for
$\theta\in\ell_0(s)$. On the other hand, if the rank of $X$ is at
least~$s$, then $R^\ast(X\ell_0(s), \sigma) = \Omega(s \sigma^2)$.
Hence, the nearest neighbor estimator (and the projected nearest
neighbor estimator) is $O(\log p)$ minimax for the hard sparsity
constraint. This is consistent with the bound in
Theorem \ref{thmmain} by letting $q\to0$.
\end{remark}

\begin{remark}
In the proof of Theorem \ref{thmlowerbound}, we actually showed that
there exists a submatrix $X'$ which consists of $k\leq n$ columns of $X$
such that the minimax risk of $X'\ell_q^k$ is close to that of
$X\ell_q^{p}$. In some sense, this means that there is a hardest
sub-problem which has at most $n$ columns.
\end{remark}

\begin{remark}
Our technique still leaves a gap of $(\log p)^{1-q/2}$.
We do not know if this gap is inherent to the projected nearest
neighbor estimator or due
to the deficiency of the analysis. We note that the upperbound cannot be
improved in general, as demonstrated by the example of $\ell_1$
ball. There might be a chance to improve the lowerbound by a factor of
$\sqrt{\log k}$ by more sophisticated techniques. But this is still
insufficient to close
the gap as $k$ might be much smaller than $p$.
\end{remark}

\begin{remark}
While PNN may sound similar to the technique of low dimension
projection, there are significant differences. For example, when
applying low dimension projection, we typically would like to preserve
the original metric structure, and often a random projection
suffices. In our case, however, we would like to make the projection
as small as possible, and it requires more careful selection of the
projection. Indeed, it is easy to show that a random projection would
fail for our purpose.
\end{remark}

\section{Algorithms}\label{secalgo}

While the analysis of projected nearest neighbor estimators is somewhat
involved, the
resulted algorithm is quite straightforward. There are two separate
parts in the projected nearest neighbor estimator. First, for given $K$
and~$\sigma$,
compute the optimal projection $P$ and $k$. Second, for any
observation $\widetilde{y}$, apply the projection and then compute the
nearest neighbor of $P^{\bot}\widetilde{y}$ to $P^\bot K$.

We will describe these two steps separately. For the first step, by
the proof of Theorem \ref{thmmain}, it suffices to compute
$d_k(K)$. This problem is however NP-hard \cite{b-02}. But since
$K=X\ell_q^p$, $\|P^\bot K\|$ must be realized at one of $p$ column
vectors of $X$ (see the proof of Lemma \ref{lemw}). Let
$V=\{x_i\dvtx i=1,\ldots,p\}$ be the $p$ column vectors of $X$. Then
computing $d_k(K)$ reduces to computing an $(n-k)$-dimensional
projection $P'$ such that $\max\{\|P'v\|\dvtx v\in V\}$ as small as
possible. This problem has been studied in \cite{vvyz-07}, and it is
shown one can compute an $O(\sqrt{\log p})$ approximation by the
semi-definite programming relaxation. The following proposition is the
main result of \cite{vvyz-07}.
%
\begin{proposition}\label{propsemi}
For any $n\times p$ matrix $X$, $0<q\leq1$, and $0\leq k\leq n$,\break
we can compute in polynomial time an $O(\sqrt{\log p})$ approximation
to $d_k(X\ell_q)$. In addition, we can compute an $(n-k)$-dimensional
subspace $P'$ in ran\-domized polynomial time such that with high
probability, $\|P'(X\ell_q)\|=\break O(\sqrt{\log p} d_k (X\ell_q))$.
\end{proposition}

As for the second step, we need to compute the nearest neighbor on
$K=X\ell_q^p$ for any given point. This can be done by convex
programming for $q=1$. Unfortunately, we do not know how to compute
it efficiently for $q<1$. So we can only claim polynomial time nearly
optimal estimator for $K=X\ell_1$, as described in
Algorithm \ref{algol1}. For description simplicity, we have described
the algorithm in which we try all $k=1,2,\ldots,n$. Since $d_k(K)$
is monotonically decreasing, the complexity can be reduced by using a
binary search. Theorem \ref{thmalgo} follows from the above
discussion.

\begin{algorithm}[t]
\caption{Nearly optimal estimator for $X\ell_1$}\label{algol1}
\begin{algorithmic}[1]
\REQUIRE design matrix $X$ and observation $\widetilde{y}$.
\ENSURE$\widehat{y}$.
\STATE Let $x_1,\ldots, x_p$ be column vectors of $X$. Denote the
set by $Y$;
\FOR{$k \in\{1,\ldots,p\}$}
\STATE Compute a projection $P_k$ such that $z_k = \|P_k Y\| = O(\sqrt
{\log p}) d_k(K)$;
\STATE Compute $r_k = k\sigma^2 + z_k \sigma\sqrt{\log p}$;
\ENDFOR
\STATE Pick $k^{\ast} = \argmin_{k} r_k$, and let $P=P_{k^\ast}$ and
$P^{\bot}$ be the subspace orthogonal to~$P$;
\STATE Compute $\widehat{y}'$ as the nearest neighbor of $P^\bot
\widetilde{y}$ to the convex hull of $\pm P^\bot x_1,\ldots, \pm
P^\bot x_p$. This can be done by using any polynomial time convex
programming algorithm.
\STATE Set $\widehat{y} = P\widetilde{y} + \widehat{y}'$.
\end{algorithmic}
\end{algorithm}

The following proof summarizes our above discussion.
\begin{pf*}{Proof of Theorem \ref{thmalgo}}
By Proposition \ref{propsemi}, we can compute an $O(\sqrt{\log p})$
approximation $d_k'$ of $d_k(X\ell_q)$. Using this approximation, we
compute
\[
R' = O\Bigl(\min_{0\leq k\leq n} \bigl(k
\sigma^2 + c_{q} d_k'^{q}
\sigma^{2-q} (\log p)^{1-q/2}\bigr)\Bigr).
\]

Since $d_k(K) \leq d_k' \leq c \sqrt{\log p} d_k(K)$ for some
constant $c>0$, we have that
\[
R_H(K,\sigma) \leq R'= c^{q}
\log^{q/2}pR_H(K,\sigma).
\]

By Theorem \ref{thmmain}, $R_H$ is an $O((\log p)^{1-q/2})$
approximation of $R^\ast$, so $R'$ is an $O((\log p)^{q/2} (\log
p)^{1-q/2})=O(\log p)$ approximation of $R^\ast$.

When $q=1$, by Proposition \ref{propsemi}, we can compute the nearly
optimal projection $P$ and use convex programming to compute the
nearest neighbor of $P\widetilde{y}$ to $PX\ell_1$. The former can be
done in randomized polynomial time and the latter in polynomial time.
\end{pf*}

\begin{remark}
The first step of the algorithm uses the semi-definite programming
relaxation to compute a nearly optimal projection of $X\ell_q$.
While it has guaranteed\vadjust{\goodbreak} approximation ratio, it can be time
consuming. In practice, the projections on the principal subspaces of
$X\ell_2$ might serve as a good heuristics.
\end{remark}

\begin{remark}
We do not have a polynomial time estimator for $0<q<1$ because of the
lack of a polynomial time algorithm for computing the nearest neighbor
to the nonconvex body of $K=X\ell_q$. While such nearest neighbor
problem is hard, for our purpose an approximate nearest neighbor is
sufficient. In addition, we only need to succeed in an average sense
as $\widetilde{y}=y+g$ for $y\in K$ and $g$ an i.i.d. Gaussian
noise. It is interesting to know if there exists an efficient
procedure in this particular setting. We note that this problem can
be formulated under the framework of the smoothed
analysis \cite{st-04}. In both cases, we are interested in minimizing
the expected performance of an algorithm (or an estimator) in the
worst case.
\end{remark}

\section{Adaptive estimator when $C$ is not given}\label{secadaptive}

The projected nearest neighbor estimator in the last section is nearly
minimax optimal once the sparsity radius is given. In this section, we
extend the same idea to design an adaptive estimator to deal with the
case when the sparsity radius is not known. Write $C=\|\theta\|_1$.
Ideally, one would like to achieve some kind of oracle inequality with
the error bound proportional to $R^\ast(X\ell_1(C),\sigma)$, that
is, the
nearly optimal risk bound assuming $C$ is available. We can only
partially achieve this goal with an extra additive term of
$\sqrt{n\log n} \sigma^2$. Here we will focus on the case of $q=1$
for the simplicity of the exposition.

Again let $K=X\ell_1$. Intuitively, the adaptive estimator will search
for the unknown $C$ at some discrete values. In view of the
upper bound in Corollary~\ref{corupperbound}, we will only try those
$C$'s which equalize the two error terms in (\ref{equpperbound}).

Define $C_k = k \sigma/ d_k(K)$ for $k=0,1,\ldots,
n/2$. $C_k$ has the following properties:
\begin{longlist}[(2)]
\item[(1)] $C_0 \leq C_1 \leq C_2 \leq\cdots$ is monotonically increasing,
since $d_k$ is nonincreasing.
\item[(2)] There is a constant $c>0$, for $C\geq C_k$,
%
\begin{equation}
\label{eqck} R^\ast\bigl(X\ell_1(C), \sigma\bigr)\geq c k
\sigma^2.
\end{equation}
This follows from Theorem \ref{thmlowerbound}.
\end{longlist}

Further, we define $P_k$ to be the $(n-k)$-dimensional projection that
realizes $d_k(K)$, that is, minimizes $\max_{1\leq i\leq n} \| Px_i\|$
among all the $(n-k)$-dimensional projection. The adaptive estimator
will estimate $\widetilde{y}_k = P_k \widetilde{y}$ against $P_k
X\ell_1(C_k)$ using the nearest neighbor estimator, starting from
$k=0$. Suppose that the outcome is~$\widehat{y}_k$. It is easy to
show that among the $n$ estimations $\widehat{y}_k$ for
$k=0,\ldots,n$, there is one that satisfies the true oracle risk
bound, that is, with high probability, there exists $0\leq k\leq n$
such that
\[
\|\widehat{y}_k - y \|^2 = O\bigl(\sqrt{\log
p}R^\ast\bigl(X\ell_1\bigl(\|\theta_1\|\bigr),\sigma
\bigr)\bigr).\vadjust{\goodbreak}
\]

Unfortunately, we cannot determine reliably which one it is. Instead,
we can only choose one which is within $O(\sqrt{n\log n}\sigma^2)$
error. This is by finding the minimum $k$ such that
$\|\widetilde{y}_k-\widehat{y}_k\|^2$ is not too large (defined
precisely later). Algorithm~\ref{algoadaptive} contains a formal
description.

\begin{algorithm}[t]
\caption{Adaptive projected nearest neighbor estimator}\label{algoadaptive}
\begin{algorithmic}[1]
\REQUIRE design matrix $X$ and observation $\widetilde{y}$.
\ENSURE estimation $\widehat{y}$.
\FOR{$k\in\{0, 1,\ldots,n/2\}$}
\STATE Compute the $(n-k)$-dimensional projection $P_k$ that
approximately minimizes $\|PX\|$;
\STATE Compute $\widetilde{y}_k = P_k \widetilde{y}$, $X_k = P_k X$,
and $\Delta_k = \max_{i} P_k x_i$;
\STATE Set $C_k = k \sigma/ \Delta_k$
\STATE Compute $\widehat{y}_k$ to be the nearest neighbor of
$\widetilde{y}_k$ on $X_k\ell_1(C_k)$
\IF{$\|\widehat{y}_k - \widetilde{y}_k\|^2 \leq(n-k)\sigma^2 +
2\sqrt{n\log n}\sigma^2$}
\STATE Set $\widehat{y} = \widehat{y}_k + P_k^\bot\widetilde{y}$
and return;
\ENDIF
\ENDFOR
\STATE Set $\widehat{y} = \widetilde{y}$.
\end{algorithmic}
\end{algorithm}

Now we will show that the estimator given in
Algorithm \ref{algoadaptive} satisfies the bound stated in
Theorem \ref{thmadaptive}. The proof requires some properties on
$\|\widehat{y}_k-\widetilde{y}_k\|^2$ as described in
Lemma \ref{lemada}. Denote by $y_k=P_k y$ and $K_k = P_k X \ell_1
(C_k)$. Let $\delta_k$ denote the $\ell_2$ distance between $y_k$ and
$K_k$, that is, $\delta_k = \min_{z\in K_k} \| y_k - z
\|$.

\begin{lemma}\label{lemada}
There are constants $c_1,c_2>0$ such that the following holds with high
probability:
\begin{longlist}[(2)]
\item[(1)] If $y_k \in K_k$, then
\[
\|\widehat{y}_k - \widetilde{y}_k\|^2
\leq(n-k)\sigma^2 + 2\sqrt{n\log n}\sigma^2.
\]

\item[(2)] If $\delta_k^2 \geq c_1 (\sqrt{n\log n} \sigma^2 + k\sigma^2
\log p)$, then
\[
\|\widehat{y}_k - \widetilde{y}_k\|^2
\geq(n-k)\sigma^2 + 2\sqrt{n\log n}\sigma^2.
\]

\item[(3)] If $\delta_k^2 \leq c_1 (\sqrt{n\log n} \sigma^2 + k\sigma^2
\log p)$, then
\[
\|\widehat{y}_k - y_k\|^2 \leq
c_2\bigl(\sqrt{n\log n}\sigma^2 + k \sigma^2
\log p\bigr).
\]
\end{longlist}
\end{lemma}

By Lemma \ref{lemada}(1) and (2), step 6 in
Algorithm \ref{algoadaptive} serves as a test for whether $y_k$ is
sufficiently separated from $K_k$. When $y_k \in K_k$, then the test
is true with high probability, and the algorithm outputs $\widehat{y}$
and returns. But when the separation between $y_k$ and $K_k$ is large
enough [$c_1(\sqrt{n\log n}\sigma^2 + k\sigma^2\log p)$], then step 6
would test false with high probability. Theorem \ref{thmadaptive} follows
from Lemma \ref{lemada}.

\begin{pf*}{Proof of Theorem \ref{thmadaptive}}
If the test at step 6 outputs false for some $k$, then by
Lemma \ref{lemada}(1), $y_k \notin K_k$. Thus $y \notin X\ell_1(C_k)$,
that is, $C \geq C_k$. By (\ref{eqck}), we have that
$R^\ast(X\ell_1(C),\sigma)\geq c k\sigma^2$.

On the other hand, if step 6 tests true for $k$, then by
Lemma \ref{lemada}(2), $\delta_k^2 \leq c_1 (\sqrt{n\log n}
\sigma^2 +
k\sigma^2\log p)$, and by Lemma \ref{lemada}(3), $\widehat{y}$
returned at step 7 satisfies that
\[
\|\widehat{y} - y\|^2 = \|\widehat{y}_k -
y_k \|^2 + k\sigma^2\leq c_2
\bigl(\sqrt{n\log n}\sigma^2 + k \sigma^2 \log p\bigr) + k
\sigma^2.
\]

We distinguish three outcomes of step 6.
\begin{itemize}
\item Step 6 tests true for $k=0$. In this case,
\[
\|\widehat{y} - y\|^2 \leq c_2\sqrt{n\log n}
\sigma^2.
\]

\item Step 6 test true for some $k>0$ and therefore is false for $k-1$.
In this case
\[
R^\ast\bigl(X\ell_1(C),\sigma\bigr)\geq c (k-1)
\sigma^2
\]
and
\begin{eqnarray*}
\|\widehat{y} - y\|^2 &\leq& c_2\bigl(\sqrt{n\log n}
\sigma^2 + k \sigma^2 \log p\bigr) + k\sigma^2
\\
&=& O\bigl(\sqrt{n\log n}\sigma^2 + R^\ast\bigl(X
\ell_1(C),\sigma\bigr)\log p \bigr).
\end{eqnarray*}

\item Step 6 is never true so step 10 is reached. In particular, the
test is false for $k=n/2$ and hence $R^\ast(X\ell_1(C),\sigma)\geq
c_1 (n/2-1)\sigma^2$ but then $\|\widehat{y}-y\|^2 = O(n\sigma
^2) = O(R^\ast(X\ell_1(C),\sigma))$.
\end{itemize}

In all the above cases, the bound in Theorem \ref{thmadaptive} holds.
\end{pf*}

\begin{remark}\label{rmkada}
When $R^\ast(X\ell_1(\|\theta\|_1),\sigma)\geq\sqrt{n}\sigma
^2$, the
bound (\ref{eqada}) in Theorem~\ref{thmadaptive} becomes a true
oracle risk bound (within $O(\log p)$ factor). In view of the proof of
Theorem \ref{thmmain}, this happens when $\|\theta\|_1
d_{\sqrt{n}}(X\ell_1) \geq
\sqrt{n}\sigma$, that is, when $\|\theta\|_1
\geq\sqrt{n}\sigma/d_{\sqrt{n}}(X\ell_1)$. In such case, the
risk ranges between $\sqrt{n\log n}\sigma^2$ and $n\sigma^2$. So
the bound (\ref{eqada}) is nearly optimal and nontrivial for a rather
large range of $\|\theta\|_1$.
\end{remark}

\begin{remark}
It might be possible to apply the Lasso or Dantzig selector estimators
to the projection $P_kX$ to obtain $\widehat{y}_k$ and then choose one
$\widehat{y}_k$ similar to Algorithm \ref{algoadaptive}. This would
probably result in the same bound as in (\ref{eqada}). We choose our
current exposition because Lemma \ref{lemada}(2) relies on the fact
that $\widehat{y}_k$ is the nearest neighbor to $P_k\widetilde{y}$. It
is not immediately clear whether it also holds for Lasso or Dantzig
selector.
\end{remark}

\begin{remark}
One may wonder if it is possible to get rid of
$\sqrt{n\log n}\sigma^2$ factor and obtain a pure oracle inequality
bound. If such a bound is possible, then when $C=0$, the estimator
needs to map all the observations to $0$. Since it is impossible to
distinguish $0$ and a sphere with radius $n^{1/4}\sigma$, there might
be a good reason for such an additive separation to be expected.
\end{remark}

\section{Proofs}\label{secproofs}

\subsection{\texorpdfstring{Proof of Lemma \protect\ref{lemnnbad} (bad example for the nearest neighbor estimator)}
{Proof of Lemma 7 (bad example for the nearest neighbor estimator)}}\label{subsecbadexample}

We will now construct a bad example for the nearest neighbor
estimator. While it is well known that the nearest neighbor estimator
can be nonoptimal, we could not find a definitive reference for a
large gap. In our example, we will demonstrate a large gap of
$\sqrt{n}$. Consider the ellipsoid
\[
E_n= \Biggl\{y=(y_1,\ldots,y_n)\dvtx  \sum
_{i=1}^{n-1} y_i^2 +
\frac
{y_n^2}{\sqrt{n}}\leq1\Biggr\}.
\]

Set $\sigma=1$. The orthogonal projection estimator $M(\widetilde{y})
= (0,\ldots,0,\widetilde{y}_n)$ has minimax error
%
\begin{equation}
\label{equ} M(\widetilde{y})=\sum_{i=1}^{n-1}
y_i^2 +\mean\bigl[(\widetilde{y}_n
-y_n)^2\bigr] \leq2.
\end{equation}

On the other hand, we show that the nearest neighbor estimator has error
$\Omega(\sqrt{n})$. For any $\widetilde{y}=(\widetilde{y}_1,\ldots,
\widetilde{y}_n)$, by using Lagrangian multiplier, we have that the
nearest point $\widehat{y}$ to $\widetilde{y}$ on $E_n$
satisfies that $\widetilde{y}_i = (1+\lambda) \widehat{y}_i$ for
$i=1,\ldots, n-1$
and $\widetilde{y}_n= (1+\lambda/\sqrt{n}) \widehat{y}_n$. Now,
pick $y=(0,\ldots, 0, n^{1/4})\in E_n$. Then with high probability
$\sum_{i=1}^{n-1} \widetilde{y}_i^2 = \Omega(n)$. By
\[
\sum_{i=1}^{n-1}\widetilde{y}_i^2
= (1+\lambda)^2 \sum_{i=1}^{n
-1}
\widehat{y}_n^2\leq(1+\lambda)^2,
\]
we have $\lambda= \Omega(\sqrt{n})$.
But then $\widehat{y}_n\leq c \widetilde{y}_n\leq c n^{1/4}$
for some
constant $c<1$. Thus, with high probability $\| \widehat{y} - y\| =
\Omega(n^{1/4})$. So the nearest neighbor estimator has error
$\Omega(n^{1/2})$. Since the projection estimator achieves the risk of
$O(1)$, we have constructed an example to show that the nearest neighbor
estimator can be $\Omega(\sqrt{n})$ factor larger than the optimal.

\subsection{\texorpdfstring{Proof of Proposition \protect\ref{proplp}}
{Proof of Proposition 9}}\label{subsecproplp}

It is well known that the error of the nearest neighbor estimator is
determined by the metric structure of $K$. For two bodies
$K_1,K_2\subseteq\reals^n$, define the (dyadic) entropy number $e_k(K_1,
K_2)$, for any $k\geq0$, as the minimum $\epsilon$ such that
$K_1$ can be covered by $2^k$ copies of $\epsilon K_2$. When $K_2$ is
the unit $\ell_2$ ball, we simply write it as $e_k(X_1)$.

For a random vector $g\in\G=\G^n(1)$ and any $y\in\reals^n$, let
$g_y$ denote the random variable $g\cdot y\in\reals$. The classical
Dudley bound states
that there is a constant $c>0$ such that
\[
\mean_{g\sim\G}\Bigl[\sup_{y\in K}|g_y|\Bigr]
\leq c \sum_{k=0}^\infty2^{k/2}
e_{2^k}(K).
\]

We need a slight variation of the above bound where the summation is
over $k$ above some threshold. For $\delta\geq0$, write
\begin{eqnarray*}
k(\delta) &=& \bigl\lfloor\log\bigl(\min\bigl\{k\dvtx  e_k(K)\leq\delta
\bigr\}\bigr) \bigr\rfloor,
\\
\gamma(K,\kappa) &=& \sum_{k=\kappa}^\infty2^{k/2}
e_{2^k}(K),
\\
K(\delta) &=& K\cap\ell_2^n(\delta).
\end{eqnarray*}

With the above notation, the following lemma holds.
%
\begin{lemma}\label{lemchain}
There is a constant $c>0$, for any $t>0$,
\[
\prob_{g\sim\G}\Bigl[\sup_{y\in K(\delta)}|g_y| \geq t
\gamma\bigl(K, k(\delta)\bigr)\Bigr] \leq\exp\bigl(-c t^2
2^{k(\delta)}\bigr).
\]
\end{lemma}
\begin{pf}
By the standard chaining argument \cite{t-05}. Clearly the result
holds if we replace $e_k(K)$ with any upper bound of $e_k(K)$.
\end{pf}

Now we prove Proposition \ref{proplp}. Without loss of
generality, we assume $\sigma=1$. We apply the standard technique to
bound the error of the nearest neighbor estimator by the supreme of
Gaussian processes \cite{vdg-00,rwy-11}. The starting point is the
well-known observation that for $\widehat{y} = N_K(\widetilde{y})$,
%
\begin{equation}
\label{eqlse} \|\widehat{y} - y\|^2 \leq2(\widetilde{y}-y)\cdot(
\widehat{y}-y).
\end{equation}

Since $\widehat{y},y\in K=X\ell_q^p$ and by the quasi-convexity of
$\ell_q^p$ for $0<q\leq1$, we have that $\widehat{y}-y\in c'
K$ for $c' = 2^{{1}/{q}}$. Observe that $g=\widetilde{y}-y$ is a
Gaussian random vector. We can bound $\|\widehat{y}-y\|$ through
Dudley bound over $\ell_q^p$ ball as follows.

To apply Lemma \ref{lemchain}, we need an estimate on the entropy
number of $K=X\ell_q^p$. Write $\Delta=\|K\|$. The following is a
consequence of \cite{cp-88,gl-00}. For completeness, we include the
derivation in the \hyperref[app]{Appendix}.
%
\begin{lemma}\label{lemmet}
%
\begin{equation}\qquad
e_{2k}\bigl(X\ell_q^p, \ell_2^n
\bigr) = \cases{ O(\Delta), &\quad $k\leq\log p$,
\vspace*{2pt}\cr
\displaystyle O \biggl( \biggl(f_q
\frac{\log(1+p/k)}{k} \biggr)^{1/q-1/2} \Delta\biggr), &\quad $\log p\leq k
\leq p$,
\vspace*{2pt}\cr
O \bigl(2^{-2k/p} (f_q/p)^{1/q-1/2}\Delta\bigr), &\quad
$k\geq p$,}
\end{equation}
where $f_q= O(\frac{1}{q}\ln\frac{2}{q})$ is a constant dependent
on $q$ only.
\end{lemma}

Now the crucial lemma is Lemma \ref{lemmain}.
%
\begin{lemma}\label{lemmain}
Suppose that $\Delta\leq p^{1/q}(\log p)^{1/2}$ and
$\Delta/p^{1/q-1/2}\leq\delta\leq\Delta$, for any constant $d>0$,
there exists $c(q,d)>0$, dependent on $q$ and $d$ only, such that
\begin{eqnarray*}
&&\prob_{g\sim\G}\Bigl[\sup_{y\in K,\|y\|\leq\delta} |g_y| \geq
c(q,d) \Delta^{{q}/({2-q})}\delta^{({2-2q})/({2-q})}\sqrt{\log p}\Bigr]
\\
&&\qquad\leq p^{-d ({\Delta}/{\delta} )^{q/2}}.
\end{eqnarray*}
\end{lemma}
\begin{pf}
The proof is by applying Lemmas \ref{lemchain} and \ref{lemmet}. By
Lemma \ref{lemmet}, for
\[
\Delta/p^{1/q-1/2}\leq\delta\leq\Delta,
\]
we have,
\[
k(\delta) = O\bigl((\Delta/\delta)^{{2q}/({2-q})} \log p\bigr) = O(p).
\]

Therefore,
\begin{eqnarray*}
\gamma\bigl(K,k(\delta)\bigr) & = &\sum_{k=k(\delta)}^\infty2^{k/2}
e_{2^k}(K)
\\
& = &\sum_{k=k(\delta)}^{\log p} 2^{k/2}
e_{2^k}(K) + \sum_{k=\log
p}^\infty2^{k/2}
e_{2^k}(K).
\end{eqnarray*}

By Lemma \ref{lemmet}, it is easily seen that for both terms, the
dominant term is the first term, that is, when $k=k(\delta)$ and
$k=\log
p$, respectively. Plugging in $e_k(K)$ for these values, we have
\begin{eqnarray*}
\gamma\bigl(K,k(\delta)\bigr) & \leq &O \bigl(\delta\sqrt{(\Delta/\delta
)^{{2q}/({2-q})}\log p} \bigr) + O \bigl(\sqrt{p}\Delta/p ^{1/q-1/2}
\bigr)
\\
& \leq &O \bigl(\Delta^{{q}/({2-q})} \delta^{
({2-2q})/({2-q})}\sqrt{\log p} +
p^{1-1/q}\Delta\bigr).
\end{eqnarray*}

It is easy to verify that with $\delta\geq\Delta/p^{1/q-1/2}$,
\[
\Delta^{{q}/({2-q})} \delta^{({2-2q})/({2-q})}\sqrt{\log p} \geq c
\Delta
p^{1-1/q} \sqrt{\log p}
\]
for some constant $c'>0$. So the first term dominates, that is,
\[
\gamma\bigl(K,k(\delta)\bigr)=O\bigl(\Delta^{{q}/({2-q})}\delta^{
({2-2q})/({2-q})}
\sqrt{\log p}\bigr).
\]

The claim now follows from Lemma \ref{lemchain}.
\end{pf}

With the above preparation, we are ready to prove Proposition \ref{proplp}.
\begin{pf*}{Proof of Proposition \ref{proplp}}
We assume $\sigma=1$. Recall $\Delta= \|K\|$. We can further assume
%
\begin{equation}
\label{eqbound} \sqrt{\log p} \leq\Delta\leq n^{1/q}/(\log
p)^{(2-q)/2q}.
\end{equation}
Otherwise the claim follows immediately by using the trivial bound
of\break
$O(\min(\Delta^2,n\sigma^2))$. Together with the assumption that
$p
=
\Omega(n/\log n)$, the upper bound in (\ref{eqbound}) implies that
%
\begin{equation}
\label{eqbound2} \Delta=O\bigl(p^{1/q}(\log p)^{1/2}\bigr).
\end{equation}

Write\vspace*{1pt} $\delta_0 = c \Delta^{q/2} (\log p)^{1/2 - p/4}$ for some
sufficiently large $c$ such that $\delta_0 \geq
\Delta/p^{1/q-1/2}$.
This is possible as $\Delta= O(p^{1/q}(\log p)^{1/2})$. Hence, by
applying Lemma \ref{lemmain}, we have that
for $\delta_0\leq\delta\leq\Delta$ and any $d>0$ there exists
$c(q,d)>0$ such that
\[
\prob_{g\sim\G} \Bigl[\sup_{y\in K(\delta)} |g_y| \geq
c(q,d) \Delta^{{q}/({2-q})}\delta^{({2-2q})/({2-q})}\sqrt{\log
p}\Bigr] \leq p
^{-d ({\Delta}/{\delta} )^{q/2}}.
\]

Now denote by $\E$ the following event:
\[
\exists y \bigl(\delta_0 \leq\|y\| \leq\Delta\bigr)
\wedge\bigl(|g_y|\geq t \sqrt{\log p}\Delta^{{q}/({2-q})} \| y
\|^{({2-2q})/({2-q})} \bigr).
\]

By the peeling argument we show that we can choose $t$, dependent on
$q$ only, such that $\prob[\E]\leq p^{-4/q}$. Define
\[
\overline{K}(\delta) = K(\delta)\setminus K(\delta/2).
\]

Clearly
$\overline{K}(\delta)\subseteq K(\delta)$ and for any $y\in
\overline{K}(\delta)$, $\|y\|\geq\delta/2$. By these we have
\[
\prob\Bigl[\sup_{y\in\overline{K}(\delta)} |g_y| \geq t_{q}
\sqrt{\log p}\Delta^{{q}/({2-q})}\|y\|^{q/2}\Bigr] \leq
p^{-d (\Delta/\delta
)^{q/2}}.
\]

Hence for any $d>0$, there is $c(q,d)>0$ such that
\begin{eqnarray*}
\prob[\E]
&=&\prob\Bigl[\sup_{y\in K, \|y\|\geq\delta_0} |g_y|\geq c(q,d) \sqrt{
\log p}\Delta^{{q}/({2-q})}\|y\|^{({2-2q})/({2-q})}\Bigr]
\\
&\leq&\sum_{k=0}^{\log(\Delta/\delta_0)} \prob\Bigl[\sup
_{y\in
\overline{K}(2^k\delta_0)} |g_y| \geq c(q,d) \sqrt{\log p}\Delta
^{{q}/({2-q})}\|y\|^{({2-2q})/({2-q})}\Bigr]
\\
&\leq&\sum_{k=0}^{\log(\Delta/\delta_0)} p^{-d(\Delta
/(2^k\delta_0))^{q/2}}.
\end{eqnarray*}

Now choosing $d=4/q$ and setting $t_{q}=c(p,4/q)$, we have that $\prob
[\E]=O(p^{-4/q})$. Let $z =
\widehat{y}-y$. So for $\|z\|\geq\delta_0$, with probability
$1-O(p^{-4/q})$,
\[
\|z\|^2\leq2|w\cdot z| \leq t_{q} \sqrt{\log p}
\Delta^{
{q}/({2-q})}\|z\|^{({2-2q})/({2-q})}.
\]

That is,
\[
\|z\| = O\bigl(\Delta^{q/2}(\log p)^{1/2-q/4}\bigr)=O(
\delta_0).
\]

Hence with probability $1-O(p^{-4/q})$,
\[
\|\widehat{y}-y\|^2 = O\bigl(\delta_0^2
\bigr) = O\bigl(\Delta^{q} (\log p)^{1-q/2}\bigr).
\]

Since $\|\widehat{y}-y\|\leq2\Delta\leq2p^{1/q}$, we have that
\[
\mean\bigl[\|\widehat{y}-y\|^2\bigr] \leq\delta_0^2
+ O\bigl(p^{-4/q}\cdot2p ^{2/q}\bigr) = O\bigl(
\Delta^{q}(\log p)^{1-q/2}\bigr).
\]

For general $\sigma>0$, we apply the standard scaling formula of
$R_N(K,\sigma) = \sigma^2 R_N(K/\sigma, 1)$ and complete the proof
of Proposition \ref{proplp}. The constant of $c_q= O(2^{
{1}/{q}}\frac{1}{q}\ln\frac{2}{q})$ comes from multiplying $c'$ and
$f_q$ in Lemma \ref{lemmet}.
\end{pf*}

%
%
%
%

\subsection{\texorpdfstring{Proof of Theorem \protect\ref{thmlowerbound}}{Proof of Theorem 11}}\label{subsecthmlowerbound}

To establish the lower bound, we consider the largest Euclidean ball
of various dimension contained in $K$. Intuitively, we show that if
Kolmogorov width of $K$ is large then it has to contain a large enough
Euclidean ball, in terms of both radius and the dimension, which allows
us to nearly match the upper bound. The crucial technical tool is
the restricted invertibility result by Bourgain and Tzafriri \cite
{bt-87} and
developed by Szarek and Talagrand \cite{st-89} and
Giannopoulous \cite{g-95}.

\begin{definition}
For a set of vectors $S$, let $\spn[S]$ denote the linear subspace
spanned by $S$. A set $V=\{v_1,\ldots, v_s\}$ is called $\delta$-wide
if for any $1\leq i\leq s$, $\dist(v_i,\spn[V/\{v_i\}])\geq\delta$,
where $\dist(v, P)$ denotes the minimum distance between $v$ and any
vector in $P$.
\end{definition}

The following proposition can be gleaned from work
in \cite{bt-87,st-89,g-95}. See \cite{jz-12} (Proposition 5.2) for a proof.
%
\begin{proposition}\label{prow}
For any $\delta$-wide set $V =
\{v_1,\ldots,v_s\}$, there exists $S \subseteq
\{1,\ldots, s\}$ with $|S|\geq(1-\epsilon) s$ such that for
any $\alpha= (\alpha_j)_{j \in S}$, $\|\sum_{j\in S} \alpha_j
v_j \| \geq c \sqrt{\epsilon/s} \delta\sum_{j\in S}
|\alpha_j|$, where $c$ is an absolute constant.
\end{proposition}

We make the following observation.
%
\begin{lemma}\label{lemw}
Suppose that $K=X\ell_q^p$ and $X=(x_1,\ldots,x_p)$.
Then for any \mbox{$k>0$}, there exists $k+1$ vectors $V\subseteq\{x_1,\ldots, x_p\}$ such that $V$ is $d_k(K)$ wide.
\end{lemma}
\begin{pf}
For a set of points $p_1,\ldots, p_s$ and $k\geq s-1$, let
$\vol_k(p_1,\ldots, p_s)$ denote the $k$-volume of the convex hull of
$p_1,\ldots, p_s$.

We find $k+1$ points $V=\{v_1,\ldots, v_{k+1}\}$ in $K$ such that the
$k+1$ volume of the simplex spanned by the origin $O$ and $v_1,\ldots,
v_{k+1}$ is the maximum, that is,
\[
V = \argmax_{y_1,\ldots, y_{k+1}\in K} \vol_{k+1}(O, y_1,\ldots,
y_{k+1}).
\]

Since $K$ is a compact set, $V\subseteq K$. We first show that $V$ is
$d_k(K)$ wide. Consider the $k$-dimensional subspace $P$ spanned by
$v_1,\ldots, v_k$. By the definition of $d_k$, we have $\sup_{y\in K}
\|Py - y\| \geq d_k(K)$. Or equivalently
%
\begin{equation}
\label{eqmax} \sup_{y\in K}\dist\bigl(y,\spn\bigl[\{v_1,\ldots, v_k\}\bigr]\bigr)\geq d_k(K).
\end{equation}

On the other hand,
%
\begin{eqnarray}\label{eqvol}
&&\vol_{k+1}(O, v_1,\ldots, v_{k+1})
\nonumber\\[-8pt]\\[-8pt]
&&\qquad= \frac{1}{k+1} \vol_{k}(O, v_1,\ldots,
v_k) \cdot\dist\bigl(v_{k+1}, \spn\bigl[
\{v_1,\ldots, v_k\}\bigr]\bigr).\nonumber
\end{eqnarray}

By the maximality of $\vol_{k+1}(O,v_1,\ldots,v_{k+1})$ and
(\ref{eqmax}) and (\ref{eqvol}), we have
\[
\dist\bigl(v_{k+1},\spn(v_1,\ldots,v_k)\bigr)
\geq d_k(K).
\]

Repeating this argument for each $v_i$ in $V$, we have that $V$ is
$d_k(K)$-wide. In addition, for $K=X\ell_1$, $K$ is the convex hull of
$\pm x_1,\ldots, \pm x_p$. Hence for any projection $P$,
$\argmax_{x\in K} \| Px\|$ has to be a vertex of $K$. That is,
$V\subseteq\{\pm x_i\dvtx  1\leq i\leq p\}$. It is easy to see that $V$
can be chosen such that $V\subseteq\{x_1,\ldots, x_p\}$. Since
$X\ell_q\subseteq X\ell_1$ for $0<q<1$, $d_k(X\ell_q)\leq
d_k(X\ell_1)$. This holds for any $0<q\leq1$.
\end{pf}

Using Proposition \ref{prow} and Lemma \ref{lemw}, we have Lemma \ref{lemradius}.
%
\begin{lemma}\label{lemradius}
There exists a constant $c>0$ such that for any $K=X\ell_q^p$, $k>0$,
and $0<\epsilon<1$, there exists a linear sub-space $P$ such that
$P\cap K$ contains an $(1-\epsilon)k$-dimensional $\ell_2$ ball with
radius $\Omega(\sqrt{\epsilon(1-\epsilon)} k^{1/2-1/q} d_k(K))$.
\end{lemma}
\begin{pf}
Clearly we can assume that $d_k(K)>0$. Let $V$ be the $d_k(K)$-wide set
as in Lemma \ref{lemw}. Write
$S_0=\{i\dvtx  x_i \in V\}$. By Proposition \ref{prow}, let $S\subseteq
S_0$ be such that $|S|\geq(1-\epsilon)|S_0|$ and for
any $\{\alpha_j\}_{j\in S}$,
\[
\biggl\| \sum_{i\in S} \alpha_i
x_i \biggr\| \geq c\sqrt{\epsilon/|S_0|} d_k(K)
\sum_{i\in S} |\alpha_i|.
\]

According\vspace*{1pt} to reverse H\"{o}lder inequality, for $x\in\reals^p$ and
$0<q\leq1$, $\|x\|_1 \geq n^{1-1/q}\|x\|_{q}$. Hence, for any $\{
\alpha_i\}$ such that $\sum_{i\in S}|\alpha_i|^{q} = 1$, $\sum_{i\in S}
|\alpha_i| \geq|S|^{1-1/q}$. Thus if $\|\alpha\|_{q} = 1$, then
%
\begin{eqnarray}\label{eqball}
\biggl\| \sum_{i\in S} \alpha_i
x_i \biggr\| &\geq& c\sqrt{\epsilon/|S_0|} d_k(K)
|S|^{1-1/q}
\nonumber\\[-8pt]\\[-8pt]
&\geq&c \sqrt{\epsilon(1-\epsilon)} |S|^{1/2-1/q} d_k(K).\nonumber
\end{eqnarray}

Let $P$ be the sub-space spanned by $x_i$ for $i\in S$. Since
$\{x_i\}_{i\in S_0}$ is $d_k(K)>0$ wide, they are linearly independent.
That is, $K\cap P$ is fully ($|S|$) dimensional. On the other hand by
(\ref{eqball}) for any $v$ on the boundary of $K\cap P$, we have that
\[
\|v\| \geq c \sqrt{\epsilon(1-\epsilon)} |S|^{1/2-1/q} d_k(K).
\]

Hence, $K\cap P$
contains an $|S|$-dimensional $\ell_2$ ball with radius
\[
c\sqrt{\epsilon(1-\epsilon)} |S|^{1/2-1/q} d_k(K).
\]
The claim follows by $|S| \leq k$ and $1/2-1/q<0$.
\end{pf}

By Lemma \ref{lemlower}, $R^{\ast}(\ell_2^k(r), \sigma) =
\Omega(\min(k\sigma^2, r^2))$. In addition, by definition of minimax
risk, for any $K_1\supseteq K_2$, $R^\ast(K_1,\sigma)\geq
R^\ast(K_2,\sigma)$ (see, e.g., \cite{dlm-90}).
Choosing $\epsilon=1/2$, we have that
for $K=X\ell_q^p$,
\[
R^\ast(K,\sigma) = \Omega\Bigl(\max_k \min\bigl(k
\sigma^2, k^{1-2/q} d_k(K)^2\bigr)
\Bigr).
\]

\subsection{\texorpdfstring{Proof of Theorem \protect\ref{thmmain}}{Proof of Theorem 1}}\label{subsecthmmain}

Let
\[
k^\ast=\argmax_k \min\bigl(d_k(K),
k^{1/q} \sigma\bigr).
\]

When there is a
tie, we pick $k^\ast$ to be the smallest among the ties. Clearly
$0<k^\ast<n$ since $d_n(K)=0$. When $k^\ast=1$, it is easy to
show the
claim holds. For $1<k^\ast<n$, we distinguish two cases.

\textit{Case} 1. $d_{k^\ast}(K)\geq(k^\ast)^{1/q}\sigma$.

In this case, we have that $d_{k^\ast+1}(K)\leq(k^\ast
+1)^{1/q}\sigma$. Otherwise, we would have that
\begin{eqnarray*}
&&\min\bigl(d_{k^\ast+1}(K), \bigl(k^\ast+1\bigr)^{1/q}
\sigma\bigr)
\\
&&\qquad= \bigl(k^\ast+1\bigr)^{1/q}\sigma> \bigl(k^\ast
\bigr)^{1/q}\sigma
\\
&&\qquad\geq d_{k^\ast}(K) \geq\min\bigl(d_{k^\ast}(K),
\bigl(k^\ast\bigr)^{1/q}\sigma\bigr).
\end{eqnarray*}

This contradicts with the maximality of $k^\ast$. Since $d_{k^\ast
}(K)\geq(k^\ast)^{1/q}\sigma$, $k^{1-2/q}\*d_{k^\ast}(K)^2 \geq
k^\ast\sigma^2$.
We apply the lower bound in (\ref{eqlowerbound}) and obtain that
\[
R^\ast(K,\sigma) = \Omega\bigl(k^\ast\sigma^2
\bigr).
\]

For the upper bound, by taking $k=k^\ast+1$ in (\ref{equpperbound}),
we have
\begin{eqnarray*}
R_H(K,\sigma)
&=&O\bigl(\bigl(k^\ast+1\bigr)\sigma^2 + c_{q}
d_{k^\ast+1}(K)^{q} \sigma^{2-q} (\log p)^{1-q/2}
\bigr)
\\
&=&O\bigl(\bigl(k^\ast+1\bigr)\sigma^2 + c_{q}
\bigl(\bigl(k^\ast+1\bigr)^{1/q}\sigma\bigr)^{q}
\sigma^{2-q}(\log p)^{1-q/2}\bigr)
\\
&=&O\bigl(\bigl(k^\ast+1\bigr)\sigma^2 (\log
p)^{1-q/2}\bigr)
\\
&=&O\bigl(R^\ast(K,\sigma) (\log p)^{1-q/2}\bigr).
\end{eqnarray*}

\textit{Case} 2. $d_{k^\ast}(K) < (k^\ast)^{1/q}\sigma$.

In this case, $d_{k^\ast}(K)\geq(k^\ast-1)^{1/q}\sigma$. Otherwise,
we would have that\break $d_{k^\ast}(K) < (k^\ast-1)^{1/q}\sigma$ and
$d_{k^\ast}(K) < d_{k^\ast-1}(K)$. The latter is due to that we pick
$k^\ast$ the smallest $k$ in case there is a tie. This would imply that
\begin{eqnarray*}
&&\min\bigl(d_{k^\ast-1}(K), \bigl(k^\ast-1\bigr)^{1/q}
\sigma\bigr)
\\
&&\qquad>d_{k^\ast}(K) \geq\min\bigl(d_{k^\ast}(K), \bigl(k^\ast
\bigr)^{1/q}\sigma\bigr).
\end{eqnarray*}

Again it contradicts with the maximality of $k^\ast$.
Hence for the lower bound, we have that
\begin{eqnarray*}
R^\ast(K,\sigma) & = & \Omega\bigl(\bigl(k^\ast
\bigr)^{1-2/q}d_{k^\ast}(K)^2\bigr)
\\
& = &\Omega\bigl(\bigl(k^\ast\bigr)^{1-2/q}\bigl(k^\ast-1
\bigr)^{2/q}\sigma^2\bigr)
\\
& = &\Omega\bigl(\bigl(k^\ast\bigr)^2 \sigma^2
\bigr)\qquad \mbox{by $k^\ast>1$.}
\end{eqnarray*}

Setting $k=k^\ast$ in (\ref{equpperbound}), we have
\begin{eqnarray*}
R_H(K,\sigma)
&=&O\bigl(k^\ast\sigma^2 + c_{q}
d_{k^\ast}(K)^{q} \sigma^{2-q} (\log p)^{1-q/2}
\bigr)
\\
&=&O\bigl(k^\ast\sigma^2 + c_{q} \bigl(
\bigl(k^\ast\bigr)^{1/q}\sigma\bigr)^{q} \sigma
^{2-q}(\log p)^{1-q/2}\bigr)
\\
&=&O\bigl(k^\ast\sigma^2 (\log p)^{1-q/2}\bigr)
\\
&=&O\bigl(R^\ast(K,\sigma) (\log p)^{1-q/2}\bigr).
\end{eqnarray*}

Therefore, for any $0<q\leq1$ and $p=\Omega(n/\log n)$, for
$K=X\ell_q^p$ where $X$ is an $n\times p$ matrix, we have that
$R_H(K,\sigma) = O((\log p)^{1-q/2} R^\ast(K,\sigma))$.

\subsection{\texorpdfstring{Proof of Lemma \protect\ref{lemada}}{Proof of Lemma 13}}

In what follows, all the statements hold with high probability, say
$1-1/n^2$.
\begin{longlist}[(2)]
\item[(1)]
Since $\widetilde{y}_k - y_k$ is $(n-k)$-dimensional Gaussian
vector, by the property of $\chi^2$-distribution,
\[
\|\widetilde{y}_k - y_k\|^2 \leq(n-k)
\sigma^2 + 2 \sqrt{n\log n}\sigma^2.
\]

Since $\|\widetilde{y}_k - \widehat{y}_k\|\leq\|\widetilde{y}_k -
y_k\|$, the statement follows immediately.

\item[(2)]
Let $z$ denote the nearest neighbor of $y_k$ on $K_k$. So $\|z-y_k\| =
\delta_k$. Further,
%
\begin{equation}
\label{eqx} (\widehat{y}_k-z)\cdot(y_k-z)\leq0.
\end{equation}

Following the same analysis for the nearest neighbor estimator, we
have
\begin{eqnarray*}
\|\widehat{y}_k - z\|^2
&\leq& 2(\widehat{y}_k - z)\cdot(\widetilde{y}_k-z)
\\
&=&2(\widehat{y}_k - z)\cdot(\widetilde{y}_k-y_k)+2(
\widehat{y}_k - z)\cdot(y_k-z)
\\
&\leq&2(\widehat{y}_k - z)\cdot(\widetilde{y}_k-y_k)
\qquad\mbox{by (\ref{eqx})}
\\
&\leq& c_1 C_k d_k \sigma\sqrt{\log p}
\\
&=&4 c_1 k\sigma^2 \sqrt{\log p}.
\end{eqnarray*}

Hence,
\begin{eqnarray*}
\|\widehat{y}_k - \widetilde{y}_k\|^2
&\geq& \|\widehat{y}_k - z\|^2 + \|z-
\widetilde{y}_k\|^2 + 2(\widehat{y}_k-z)
\cdot(z-\widetilde{y}_k)
\\
&\geq& \|\widetilde{y}_k - z\|^2 + 2(
\widehat{y}_k-z)\cdot(z-y_k) +2(\widehat{y}_k-z)
\cdot(y_k-\widetilde{y}_k)
\\
&\geq& \|\widetilde{y}_k - z\|^2 + 2(
\widehat{y}_k-z)\cdot(y_k-\widetilde{y}_k)
\qquad\mbox{by (\ref{eqx})}
\\
&\geq& \|\widetilde{y}_k - z\|^2 - 2\bigl|(
\widehat{y}_k-z)\cdot(y_k-\widetilde{y}_k)\bigr|.
\end{eqnarray*}

We bound these two terms separately:
\begin{eqnarray*}
\|\widetilde{y}_k - z\|^2
&=&\|\widetilde{y}_k - y_k\|^2 +
\|y_k-z\|^2 + 2(\widetilde{y}_k -
y_k)\cdot(y_k-z)
\\
&\geq&(n-k)\sigma^2 - 2 \sqrt{n\log n}\sigma^2 +
\delta_k^2 - 4 \delta_k \sigma\sqrt{\log p}.
\end{eqnarray*}

By the analysis for the nearest neighbor estimator, we have
\[
2\bigl|(\widehat{y}_k-z)\cdot(y_k-\widetilde{y}_k)\bigr|
\leq c_1 C_k d_k \sigma\sqrt{\log p} =
c_1 k \sigma^2 \sqrt{\log p}.
\]

Putting them together, we can take $\delta_k^2 = c_2 (\sqrt{n\log
n}\sigma^2 + k\sigma^2 \sqrt{\log p})$ for some sufficiently
large $c_2$ and obtain
\[
\|\widehat{y}_k - \widetilde{y}_k\|^2
\geq(n-k)\sigma^2 + 2\sqrt{n\log n}\sigma^2.
\]

\item[(3)] If $\delta_k^2 \leq c_1 (\sqrt{n\log n} \sigma^2 + k\sigma^2
\log p)$, then according to the above
\[
\|\widehat{y}_k - z\|^2\leq O\bigl(k\sigma^2
{\log p}\bigr) = O\bigl(\delta_k^2\bigr).
\]

Hence,
\[
\|\widehat{y}_k - y_k \| \leq\|\widehat{y}_k
- z\|+\|z-y_k\| = O(\delta_k).
\]
\end{longlist}

\begin{appendix}\label{app}
\section*{\texorpdfstring{Appendix: The entropy number of $X\ell_q$}{Appendix: The entropy number of X l q}}

By Guedon and Litvak (\cite{gl-00}, Theorem 6)
%
\begin{equation}
e_k\bigl(\ell_q^p, \ell_1^p
\bigr) = \cases{ \Theta(1), &\quad $k\leq\log p$,
\vspace*{2pt}\cr
\displaystyle \Theta\biggl(
\biggl(f_q\frac{\log(1+p/k)}{k} \biggr)^{1/q-1} \biggr), &\quad $\log p
\leq k \leq p$,
\vspace*{2pt}\cr
\displaystyle \Theta\bigl(2^{-k/p} (f_q/p)^{1/q-1}
\bigr), &\quad $k\geq p$,}
\end{equation}
where $f_q= O(\frac{1}{q}\ln\frac{2}{q})$ is a constant dependent
on $q$ only.

And by Carl and Pajor \cite{cp-88},
%
\begin{equation}
e_k\bigl(X\ell_1^p, \ell_2^n
\bigr) = \cases{ O(\Delta), &\quad $k\leq\log p$,
\vspace*{2pt}\cr
\displaystyle O \biggl( \biggl(\frac{\log(1+p/k)}{k}
\biggr)^{1/2} \Delta\biggr), &\quad $\log p\leq k \leq p$,
\vspace*{2pt}\cr
O
\bigl(2^{-k/p} (1/p)^{1/2} \Delta\bigr), &\quad $k\geq p$.}
\end{equation}

From the definition of $e_k$, we have (see also
\cite{pisier-89})
%
\begin{equation}
\label{eqent} e_{k_1+k_2}(K_1,K_3) \leq
e_{k_1}(K_1,K_2)e_{k_2}(K_2,K_3).
\end{equation}

By (\ref{eqent}), $e_{2k}(X\ell_q^p,\ell_2^n) \leq
e_k(\ell_q^p,\ell_1^p) e_k(X\ell_1^p, \ell_2^n)$. So we have
%
\begin{equation}\qquad
e_{2k}\bigl(X\ell_q^p, \ell_2^n
\bigr) = \cases{ O(\Delta), &\quad $k\leq\log p$,
\vspace*{2pt}\cr
\displaystyle O \biggl( \biggl(f_{q}
\frac{\log(1+p/k)}{k} \biggr)^{1/q-1/2} \Delta\biggr), &\quad $\log p\leq k
\leq p$,
\vspace*{2pt}\cr
O \bigl(2^{-2k/p} (f_{q}/p)^{1/q-1/2}\Delta\bigr), &\quad $k
\geq p$.}
\end{equation}
\end{appendix}

\section*{Acknowledgments}

The author would like to thank Adel Javanmard and Tong Zhang for useful
discussions and anonymous reviewers and the Associate Editor for many
useful suggestions, especially for suggesting the extension to the
adaptive estimation and the references
of \cite{l-adaptive-91,b-lepski-01}.



\printaddresses

\end{document}